\documentclass[12pt]{article}

\usepackage{graphicx}
\usepackage{amsmath,amsthm,amssymb,color}
\usepackage[noadjust,sort]{cite}
\usepackage[mathscr]{euscript}
 \usepackage[normalem]{ulem} 
\usepackage{subcaption}

\normalsize

\theoremstyle{plain}

\theoremstyle{definition}

\def\calF{\mathcal{F}}

\def\calL{\mathcal{L}}

\def\calR{\mathcal{R}}

\theoremstyle{definition}

\newtheorem{Problem}{Problem}[section]

{\par\noindent{\it Proof of}} % команды дл  \begin
{\hfill$\vspace{5mm}\scriptstyle\blacksquare$} % команды дл  \end

\numberwithin{equation}{section} %„тобы нумераци  в каждой секции была независимой
\def\supp{\operatorname{supp}}

\def\R{\mathcal{R}}

\def\Err{\mathfrak{E}}

\def\Reals{{\mathbb{R}}}

\def\Naturals{{\mathbb{N}}}

\def\st{\,:\,}

\def\dfrac#1#2{\lower0.15ex\hbox{\large$\textstyle\frac{#1}{#2}$}}

\begin{document}

\setcounter{page}{1}

\markboth{M. Isaev, R.G. Novikov, G.V. Sabinin}{
Numerical reconstruction from the Fourier transform on the ball}

\title{
Numerical reconstruction from the Fourier transform on the ball
using prolate spheroidal wave functions
%
%\footnotemark[3] 
% \thanks{The first author's research is  supported by    the Australian Research  Council  Discovery Early Career Researcher Award DE200101045. \\
% The work was started during the stage of the third author in the Centre de Math\'ematique Appliqu\'ees of Ecole Polytechnique in August-October 2021.
%The work is supported by a joint grant of the Russian Foundation for Basic Research and CNRS (projects no. RFBR 20-51-1500/PRC no. 2795 CNRS/RFBR).
%}
}
\date{}

\author{ 
Mikhail Isaev
\thanks{Research supported by ARC  DE200101045.}\\
\small School of Mathematics\\[-0.8ex]
\small Monash University\\[-0.8ex]
\small Clayton, VIC, Australia\\
%\small Moscow Institute of Physics and Technology\\[-0.9ex]
%\small Dolgoprudny, Moscow Region, Russia\\[-0.3ex]
\small\texttt{mikhail.isaev@monash.edu}
\and
Roman G. Novikov
\thanks{Research supported by a joint grant of the Russian Foundation for Basic Research and CNRS (projects no. RFBR 20-51-1500/PRC no. 2795 CNRS/RFBR)}
\\
\small CMAP, CNRS, Ecole Polytechnique\\[-0.8ex]
\small Institut Polytechnique de Paris\\[-0.8ex]
\small Palaiseau, France\\
\small IEPT RAS, Moscow, Russia\\
\small\texttt{novikov@cmap.polytechnique.fr}
\and
Grigory V. Sabinin\footnotemark[2]\\
\small Faculty of Mechanics and Mathematics,\\[-0.8ex]
\small Lomonosov Moscow State University\\[-0.8ex]
\small Moscow, Russia\\
\small\texttt{gvsabinin@gmail.com}
}

\maketitle
%{\bf Abstract}
%\begin{abstract}
%Many important inverse problems are exponentially ill-posed in general, which constitutes a severe difficulty for numerical treatments. However, a stable reconstruction of the unknown parameter might still be possible in some cases when the parameter is well-behaved.  Countless results in the literature confirm improved stability under various additional a-priori assumptions.  In fact, the behaviour of the stability bounds can change dramatically from the logarithmic type to the H\"{o}lder type or even, under some strong assumptions, to the Lipchitz type. In this work, we illustrate such transitions with an example from the classical Fourier analysis.
%\end{abstract}

\begin{abstract}
We implement numerically formulas of [Isaev, Novikov,  arXiv:2107.07882]
for finding a compactly supported function $v$ on $\Reals^d$, $d\geq 1$, from its Fourier transform $\calF [v]$ given within the ball $B_r$.
For the one-dimensional case, these formulas are based on the theory of prolate spheroidal wave functions,
which arise, in particular, in the singular value decomposition of the aforementioned band-limited Fourier transform for $d = 1$.
In multidimensions,  these formulas also include inversion of the  Radon transform. In particular, we give numerical examples of super-resolution, that is, recovering details beyond the diffraction limit. 
\end{abstract}

\section{Introduction}\label{S:intro}

We consider the Fourier transform $\calF$ defined by the formula

\begin{equation}
\label{eq:Fourier}
\calF[v] (p) = \hat{v}(p) :=
\dfrac{1}{(2\pi)^d}\int\limits_{\mathbb{R}^d} e^{i pq } v(q) dq, \qquad p\in \mathbb{R}^d,
\end{equation}
where $v$ is a complex-valued test function on $\Reals^d$, $d\geq 1$.

For any $\rho>0$, let
\begin{equation}
B_{\rho} := \left\{q\in \mathbb{R}^d : |q| < \rho \right\}.
%\text{for any $\rho>0$.} \nonumber
\end{equation}

We consider the following problem.

\begin{Problem}\label{P1} Find $v \in \calL^2(\Reals^d)$, where $supp\,v\subset B_{\sigma}$,
from $\hat{v} = \calF[v]$ given on the ball $B_r$ (possibly with some noise), for fixed $ r,\sigma>0$.
\end{Problem}

Problem \ref{P1} arises in different areas such as Fourier analysis, linearised inverse scattering and image processing, and has been extensively studied in the literature.
 Solving Problem~\ref{P1}  is complicated considerably  by the fact that it is   ill-posed 
in \ the sense of Hadamard (for example, when the noisy  data is taken from $\calL^2(B_r)$) and, moreover, it is exponentially unstable. Nevertheless, there exist several techniques to approach this problem theoretically and numerically. For more background on Problem \ref{P1}   see, for example, 
\cite{AMS2009,BM2009,CF2014,Gerchberg1974,IN2020,IN2020+,INnotePSWF, LRC1987,Papoulis1975} and references therein.
In addition, for general background on ill-posed inverse problems  see  \cite{TA1977, HR2021}.

% For the introduction to the theory of ill-posed  problems; see the classical  books
%  by Tikhonov, Arsenin  \cite{TA1977} and by Lavrent'ev et al.  \cite{LRS1986}. 

The conventional approach  for solving Problem  \ref{P1}   is based on the following approximation
 \begin{equation}\label{eq:naive}
 	v \approx v_{\rm naive} := \calF^{-1} \left[w \right]  \text{ on } B_{\sigma},
 \end{equation}
 where $ \calF^{-1} $ is  the standard inverse Fourier transform and $w$ is such that   $w|_{B_{r}}$ coincides with the data of Problem \ref{P1}
 and  $w|_{\Reals^d \setminus B_{r}} \equiv 0$.  
 Formula \eqref{eq:naive} leads to a stable and accurate reconstruction for  sufficiently   large $r$.   However, it  has the well-known diffraction limit: small details  (especially less than  $\pi /r$) are blurred.  A new approach for  \emph{super-resolution} in comparison with the resolution of  \eqref{eq:naive} 
 was  recently suggested in  \cite{INnotePSWF}; see also 
 \cite[Section 6.3]{IN2020}.  In the present work, we   study  numerically the approach of \cite{INnotePSWF}  and %, in particular, 
  demonstrate its efficiency.  
 %for a better  reconstruction of  small details. 

For convenience, we consider the scaling of $v$ with respect to the size of  its support:
\begin{equation}\label{def_vsigma}
	v_{\sigma} (q) := v(\sigma q), \qquad q \in \Reals^d.
\end{equation}
Note that $\supp\,v_{\sigma}\subset B_1$.
Let 
\begin{equation}\label{def_c}
	 c := r\sigma.
\end{equation}
Then,  the data in Problem \ref{P1} (for the case without noise) can be presented as follows:
%\begin{equation}\label{eq:grtheta}
% \hat{v} (r x \theta ) =
%\dfrac{\sigma}{ (2\pi) ^d}
%\int_{-1}^1 e^{i c x y } \mathcal{R} [v] (\sigma y,\theta) dy, \qquad x \in [-1,1], \ \theta \in \mathbb{S}^{d-1},
%\end{equation}
\begin{equation}\label{presentation_1D}
 \hat{v} (r x) =  \dfrac{\sigma}{2\pi}\calF_c \left[  v_\sigma \right](x),  \qquad  \text{for $d=1$},
\end{equation}
\begin{equation}\label{presentation_mult}
 \hat{v} (r x \theta ) =
\left(\dfrac{\sigma}{ 2\pi}\right)^d
 \calF_c \left[ \mathcal{R}_{\theta} [v_\sigma ]\right](x),  \qquad   \text{for $d\geq 2$,}
\end{equation}
where $x \in [-1,1]$, $\theta \in \mathbb{S}^{d-1}$; see \cite[Theorem 1.1 and Section 4.1]{INnotePSWF}. 
Here, the operators  $ \calF_c$ and  $\mathcal{R}_{\theta}$ are defined by 
\begin{align}
	\calF_c[f] (x) &:= \int_{-1}^1 e^{i c xy} f(y)dy, \qquad x \in [-1,1], 	\label{def:Fc}\\
		\mathcal{R}_\theta [u] (y)  &:= \int_{q\in \Reals^d \st q \theta =y } u(q) dq, \qquad y\in \Reals,\label{def:Rtheta}
\end{align}
where $f$ is a test function on $[-1,1]$ and  $u$ is a test function on $\Reals^d$.

Recall  that  $\mathcal{R}_{\theta}[u] \equiv \mathcal{R} [u] (\cdot, \theta)$,
where  $\mathcal{R}_{\theta}$ is defined  by \eqref{def:Rtheta}   and $\mathcal{R}$   is the classical Radon transform; see,  for example, \cite{Naterrer2001} and references therein. 
In fact, presentation \eqref{presentation_mult}  follows from the projection theorem of the Radon transform theory.
% Note also that presentation \eqref{presentation} for $d=1$ reduces to 
%$
%	 \hat{v} (r x \theta ) =  \dfrac{\sigma}{ 2\pi}  \calF_c \left[  v_\sigma \right](x)
%$

The operator $\calF_c$ defined by \eqref{def:Fc} is a variant of band-limited Fourier transform.  This operator  is one of  the key objects of  the theory of  \emph{prolate spheroidal wave functions} (PSWFs); see, for example, \cite{Slepian1983,INnotePSWF, Wang2010, BK2017, XRY2001} and references therein.  In particular, the operator $\calF_c$ has the following singular value decomposition in $\calL^2([-1,1])$:
\begin{equation}\label{Fc-dec}
\calF_c [f] (x) = \sum_{j \in \Naturals} \mu_{j,c}\psi_{j,c}(x) \int_{-1}^1 \psi_{j,c} (y) f(y) dy, 
\end{equation}
where $(\psi_{j,c})_{j \in \Naturals}$ are the 
 prolate spheroidal wave functions (PSWFs) and the    eigenvalues $\{\mu_{j,c}\}_{j \in \Naturals}$ satisfy   $0<|\mu_{j+1,c}| < |\mu_{j,c}|$ for all $j \in \Naturals$.  Here and throughout the paper, we set 
 $\Naturals : = \{0,1,2\ldots\}$.  
%These formulas are based on presentations \eqref{presentation_1D}, \eqref{presentation_mult}, inversion of $\calF_{c}$, and inversion of $\calR$. 

The   approach  for  solving Problem \ref{P1} suggested in \cite{INnotePSWF} 
 is based on presentations \eqref{presentation_1D}, \eqref{presentation_mult}, inversion of $\calF_{c}$, and inversion of $\calR$. 
The inversion of $\calR$ is given using standard results of the Radon transform theory.
The inversion of $\calF_{c}$ is given using the   singular value decomposition \eqref{Fc-dec}.
In the framework of this approach,  the operator $\calF_{c}^{-1}$ is approximated by
the finite-rank operator $\calF_{n,c}^{-1}$ (see   \eqref{def:Fnc} for precise definition), where $n$ is the rank.
In fact, the number $n$ is a regularisation parameter and its choice is crucial for  both theoretical  results and numerical applications.

We test different principles for choosing $n$, including residual minimisation and  the Morozov discrepancy principle.
One of the most interesting points of our 
numerical results lies in  examples of   super-resolution,
that is,  recovering details of size less than  $\pi/r$, where $r$ is the band-limiting radius of Problem \ref{P1}.
We also obtain a better reconstruction  in the sense of $\calL^2$-norm
than the conventional reconstruction based on formula  \eqref{eq:naive}.

%
%and $(\psi_{j,c})_{j \in \Naturals}$ coincide with \emph{prolate spheroidal wave functions}
%which are are real-valued and form an orthonormal basis in $\calL^2([-1,1])$.
%See, for example, [Isaev, Novikov,  arXiv:2107.07882] and references therein.

The  paper is structured as follows. In  Section \ref{S:recall}, we recall 
the aforementioned reconstruction formulas of \cite{INnotePSWF}. In Section \ref{S:impl}, we discuss numerical principles for choosing  the regularisation parameter $n$.   Numerical  examples   are presented  in  Section \ref{S:Examples}. 
%
%Numerical results of the present work include examples of super-resolution reconstruction,
%that is, examples of recovering details beyond the diffraction limit,
%that is, details of size less than  $\pi/r$, where $r$ is the band-limiting radius of Problem \ref{P1}.
%The numerical reconstruction developed in the present work also gives better results in the sense of $\calL^2$-norm
%than the conventional reconstruction based on formula (1.3). 
%For details, see Section 3.

%
%In our numerical examples the use of $\calF_{n,c}^{-1}$ is not less efficient than approximating $\calF_{c}^{-1}$
%using the simplest Tikhonov regularisation.
%See Sections 2, 3, and 4.
%% 

\section{Reconstruction for Problem \ref{P1}} \label{S:2}

In this section, we present the main points of our numerical approach to Problem \ref{P1}.  Namely, we recall the reconstruction formulas from  \cite{INnotePSWF} in Section  \ref{S:recall} and 
suggest their possible regularisations  in Section  \ref{S:impl}.

\subsection{Reconstruction formulas from \cite{INnotePSWF}}\label{S:recall}

 Recall the definitions of $v_\sigma$ and $c$  from \eqref{def_vsigma} and  \eqref{def_c}. 
 For the case without noise, the following reconstruction formulas for  Problem \ref{P1} hold;  see  \cite[Theorem 1.1,  Remark 1.2, and formula (1.3)]{INnotePSWF}. 
 \begin{itemize}
 	\item For $d=1$, we have
 \begin{equation}\label{inverse_1D}
 	v_\sigma =  \dfrac{2\pi}{\sigma} \calF_{c}^{-1} [\hat{v}_r],
 \end{equation}
 where 
 \[\hat{v}_r(x) = \hat{v}(rx), \qquad x\in [-1,1].
 \]
 	\item For $d\geq 2 $, we have
 \begin{equation}\label{inverse_mult}
 	v_\sigma =  \left(\dfrac{2\pi}{\sigma}\right)^d  \R^{-1}[f_{r, \sigma}], 
 \end{equation}
 where  $ \R^{-1}$ is a standard inversion of the  Radon transform $\R$, and
 \begin{align*}
 		  f_{r,\sigma}(y,\theta) &:=
 	 \begin{cases}
   \calF_c^{-1}[\hat{v}_{r,\theta}](y), &\text{if } y\in [-1,1]\\
 0, & \text{otherwise},  
 \end{cases} 
 \\
 \hat{v}_{r,\theta} (x)&:=      \hat{v} (r x\theta), \qquad  x\in [-1,1],\ \theta \in \mathbb{S}^{d-1}.
 \end{align*}
 \end{itemize}
 In the above,
   the inverse transform  $\calF_{c}^{-1}$ is given by
 	\begin{equation}\label{f:inverse}
	 \calF_{c}^{-1} [g](y) =    \sum_{j \in \Naturals} \dfrac{1}{\mu_{j,c}}\psi_{j,c}(y) \int_{-1}^1 \psi_{j,c} (x) g(x)dx,
\end{equation}
where $g$ is a test function from the range of $\calF_{c}$ acting on   $\calL^2([-1,1])$.
 
 For the case of noisy data in Problem \ref{P1}, the operator
 $ \calF_c^{-1}$ is approximated by the finite rank operator  $\calF_{n,c}^{-1}$ defined by
 \begin{equation}\label{def:Fnc}   
 	\calF_{n,c}^{-1} [g] (y) :=  \sum_{j=0}^n \dfrac{1}{\mu_{j,c}}\psi_{j,c}(y) \int_{-1}^1 \psi_{j,c} (x) g(x)dx. 
 \end{equation}
The operator   $\calF_{n,c}^{-1}$ is  correctly defined  on $\calL^2([-1,1])$ for any $n \in \Naturals$.   
In addition, $\calF_{n,c}^{-1} [g]$  is   \emph{the quasi-solution} in the sense of 
Ivanov   of  the equation 
$\calF_{c} [f] = g \in \calL^2([-1,1])$ on
%the space  $\pi_{n,c} (\calL^2([-1,1])) $
%where
%  $\pi_{n,c}[\cdot]$ is the orthogonal projection  
   %onto 
   the span of   the first $n+1$ functions
 $(\psi_{j,c})_{j  \leq  n}$.  
% \begin{equation}\label{def:pi_n}
% 	\pi_{n,c}[f]:=  \sum_{j =0 }^n  \hat{f}_{j,c}  \psi_{j,c},  
%\qquad
% 	\hat{f}_{j,c} := \int_{-1}^1 \psi_{j,c} (y) f(y) dy.
% \end{equation}

The rank $n$ of   the operator $\calF_{n,c}^{-1}$ is a regularisation parameter.   
The optimal choice of $n$ depends, in particular, on 
the relative noise level $\delta$ in the data $w \approx \hat{v}|_{B_r}$
of Problem~\ref{P1}.  
In \cite{INnotePSWF},  the pure mathematical choice of $n = n^*_{\alpha,\delta}$ is as follows:  
  \begin{equation}\label{def:n-star}
   n^*_{\alpha,\delta}  :=    \left\lfloor 3+ \tau \dfrac{ec}{4}\right\rfloor, 
 \end{equation}
   where $\lfloor\cdot \rfloor$ denotes the floor function and
   $\tau = \tau(c,\alpha,\delta) \geq 1$ is the solution of  the equation
   \begin{equation*} %\label{eq-tau}
   \tau \log \tau =    \dfrac{4}{ec}  \alpha  \log  (\delta^{-1}).
  \end{equation*}
  Here,    $\delta \in (0,1)$ is defined   using 
  $\calL^2$-norm  for $d=1$ and a weighted $\calL^2$-norm for  $d\geq 2$, 
  and     $\alpha \in (0,1)$   is  a parameter in the  related stability estimate  
for the reconstruction via formulas   \eqref{def_vsigma}, \eqref{inverse_1D}, \eqref{inverse_mult}, and \eqref{def:Fnc}; see \cite[Theorem 1.4]{INnotePSWF} for details.

In the next section, 
we discuss   numerical principles for chosing the regularisation parameter $n$.

% Let 
% \begin{equation}\label{def:pi_n}
% 	\pi_{n,c}[f]:=  \sum_{j =0 }^n  \hat{f}_{j,c}  \psi_{j,c},  
%\qquad
% 	\hat{f}_{j,c} := \int_{-1}^1 \psi_{j,c} (y) f(y) dy.
% \end{equation}
% That is,  $\pi_{n,c}[\cdot]$ is the orthogonal projection  in $\calL^2([-1,1])$  onto the span of   the first $n+1$ functions
%  $(\psi_{j,c})_{j  \leq  n}$.  

\subsection{Numerical implementation}\label{S:impl}

In this section we  describe the numerical implementation 
of  formulas \eqref{def_vsigma}, \eqref{inverse_1D},  \eqref{inverse_mult}.   
We replace  $\calF_c^{-1}$ with the finite rank operator $\calF_{n,c}^{-1}$ defined by \eqref{def:Fnc}. For implementing $\calR^{-1}$, we use the filtered back projection (FBP) algorithm, see, for example, \cite[Chapter 5]{Naterrer2001}.

The key point of our reconstruction  is  choosing the regularisation parameter $n$ in 
\eqref{def:Fnc}.
%
%In particular, we  discuss  the choice of 
%regularisation parameters.  For  $\calF_c^{-1}$, we consider two approaches: 
%\begin{itemize}
%	\item   $\calF_c^{-1}$ replaced by the finite rank operator $\calF_{n,c}^{-1}$;
%	\item $\calF_c^{-1}$ replaced by  its Tikhonov regularization $\mathcal{T}_{\delta}$.
%\end{itemize}
%For $\calR^{-1}$, we use the filtered back projection (FBP) algorithm, see, for example, \cite[Chapter 5]{Naterrer2001}.
%\[
%	\mathcal{T}_{\delta} [w] :=  
%		(\calF_c^* \calF_c +\delta \mathcal{I})^{-1}[\calF_c^*[w]] =
%	  \sum_{j \in \Naturals} \dfrac{\bar{\mu}_{j,c}}{|\mu_{j,c}|^2 + \delta}\,\psi_{j,c}(y) \int_{-1}^1 \psi_{j,c} (x) w(x)dx,
%\]
%where 
%$\mathcal{I}$ is the identity operator, $\calF_c^*$ is the  Hermitian conjugate of $\calF_c$, $\bar{\mu}_{j,c}$ is the complex conjugate of $ \mu_{j,c}$.  Equivalently,  one can define  $\mathcal{T}_{\delta} [w]$ as the function $f$ minimising   the following functional:
%\[
%	\mathcal{T}_{\delta} [w]  =
%	 \operatorname{arg} \min_{f\in \calL^2([-1,1])}\left( \|\calF_c[f] - w\|_{\calL^2[(-1+1)]}^2 + \delta \|f\|_{\calL^2([-1,1])}^2\right). 
%\]
The first interesting option  is  $n = n_0$,  where 
\begin{equation}\label{def:n0}
	n_0:= \left\lfloor \dfrac{2 c}{\pi}\right\rfloor.
\end{equation}
%where  $\lfloor \cdot \rfloor$ is the floor function. 
This choice is motivated by the following well-known formula; see, for example, \cite[formulas (2.3) and (2.4)]{INnotePSWF}:
\begin{equation}\label{eigen:ineq}
  \left\lfloor\frac{2c}{\pi}\right\rfloor-1\leq \Big|\{n \in \Naturals \st  |\mu_{n,c}|\geq \sqrt{\pi/c}\}\Big|\leq  \left\lceil \frac{2c}{\pi}\right\rceil+1.
 \end{equation}	
 In the above, $\lfloor \cdot \rfloor$ and  $\lceil \cdot \rceil$  denote the floor and the ceiling functions, respectively,  and $|\cdot|$ is the number of elements in a set.
 In fact, $|\mu_{n,c}|$ gets very small soon after $n$ exceeds $n_0$  and further 
 decays  super-geometrically as $n$ grows; see, for example, \cite{Wang2010, BK2014, BK2017, KRD2021}.
  In addition, for all our numerical examples, we observed that  
 $\calF_{n,c}^{-1}$ with $n = n_0$ leads to a reconstruction that behaves similarly to \eqref{eq:naive}.
% \[
% 	\calF_c^{-1}[w] (x) \approx coeff* \int_{-1}^{1} e^{icxy} w(y)dy
% \]
% \[
% 		\calF_c^{-1}[w] (x)  \approx coeff*\sum_{j=-n_0^*}^{n_0^*}  e^{   i c j \pi x } w(j/n_0^*)
% \]

The choice $n = n_0$  can be also intuitively explained using approximation of  $\calF_c^{-1}$ with the inverse Fourier series.  Indeed, we have that,  for $x\in [-1,1]$,  
\begin{align*}
	   \calF_c^{-1}[g] (x)&=  \pi \sum_{k= -\infty}^{\infty}   e^{ - i \pi k x}   \hat{f} (\pi k)
	\\ &\approx \pi  \sum_{k \st \pi k \in [-c,c]}  e^{ - i \pi k x}   \hat{f} (\pi k)
	 =    \dfrac{1}{2}\sum_{k \st \pi k \in [-c,c]}  e^{ - i \pi k x}  g (\pi k/c),
\end{align*}
where $g = \calF_c[f]$ and  the truncation of the series corresponds  to the known values  $ \{\hat{f} (\pi k)\}$ from $g$ given on $[-1,1]$.  Observe that $n_0$    almost  coincides with the number of terms in the truncated series above,  that is,
 the number of harmonics of the form $e^{-i \pi kx}$ periodic with respect to $x\in [-1,1]$  such that $ \pi k \in [-c,c]$.

Note that  our implementations  rely on approximations  $\{\tilde{\mu}_{j,c}\}_{j \in \Naturals}$
 of the eigenvalues $\{\mu_{j,c}\}_{j \in \Naturals}$,  approximations   $\{\tilde{\psi}_{j,c}\}_{j \in \Naturals}$  of the PSWFs $\{\psi_{j,c}\}_{j \in \Naturals}$, and methods of computing  integrals (the numbers of grid points, for example).
For consistency, we use  the tilde notation for numerical implementation of all objects and operators of our reconstruction; for example, $\tilde{\calF}_{n,c}^{-1}$ and $\tilde{\calR}^{-1}$ correspond to 
$\calF_{n,c}^{-1}$ and $\calR^{-1}$, respectively.

The quality of numerical implementations  restricts how large  $n$ could be. 
 A very rapid decay of    $  |\mu_{j,c}|$    for large $j > n_0$ 
  leads to  that dividing  by $\mu_{j,c}$ in \eqref{def:Fnc}   quickly becomes numerically intractable.
 Our \emph{trust criteria}
is 
%$n \leq  \tilde{n}_{\varepsilon}$, where $\tilde{n}_{\varepsilon}$  is  the maximum $n$ such that  $\varepsilon_n \leq \varepsilon$ 
\begin{equation}\label{Trust}
	n \leq  \tilde{n}_{\varepsilon}: = \max \{j \in \Naturals \st \varepsilon_j \leq \varepsilon\},
\end{equation}
where
$\varepsilon$ is a fixed small positive number and
\begin{equation}\label{def:epsn}
	\varepsilon_j:= \left(\sum_{\ell = 0}^j    \int_{-1}^1  \left|  \frac{\tilde{\calF}_c [\tilde{\psi}_{\ell,c}](x)}{\tilde{\mu}_{\ell,c}} - \tilde{\psi}_{\ell,c}(x)\right|^2  dx \right)^{1/2}  
\end{equation}  The  integration  and arithmetic operations  in \eqref{def:epsn} are considered as their numerical realisations. 
%The consideration of $n >\tilde{n}_{\varepsilon}$   leads to unreliable reconstructions. 

In order to  choose $n$ optimally within the window
 \begin{equation}\label{def:window}
	n_0\leq n \leq \tilde{n}_{\varepsilon}, 
\end{equation}
we rely on  the following two  well-known numerical principles.
Let $ \tilde{\calF}$ denote  the numerical implementation of the Fourier transform $\calF$ (as  explained above)  and $\tilde{\varPhi}_n[w]$  denote the numerical reconstruction 
 for Problem~\ref{P1} from the  data $w \approx \hat{v}|_{B_r}$  
 via formulas   \eqref{def_vsigma},  \eqref{inverse_1D}, \eqref{inverse_mult}, and \eqref{def:Fnc}. 
\emph{The  residual minimisation  principle} suggests 
  $ 	n = n^* $, where
  \begin{equation}\label{def:rmp}
 	 n^* := \operatorname{arg} \min_{n_0 \leq n \leq \tilde{n}_{\varepsilon}} 
 	\|
 	  \tilde{\calF}  [\tilde{ \varPhi}_n[w]]  - w
 	\|_{\calL^2(B_r)}.
 \end{equation}
 \emph{The Morozov discrepancy principle}
 suggests    $n = n^*_{\varDelta}$,  where
\begin{equation}\label{def:Morozov}
		n^*_{\varDelta} :=   
		\operatorname{arg} \min_{n_0 \leq n \leq \tilde{n}_{\varepsilon}} 
	\left|
 	\|
 	  \tilde{\calF}  [ \tilde{\varPhi}_n[w]]  -w
 	\|_{\calL^2(B_r)} -   \varDelta  \right|
\end{equation}
and
$\varDelta>0$ is a priori bound  on the  $\calL^2$ noise level of the data $w$  in Problem~\ref{P1}, that is,
\begin{equation}\label{def:delta}
	\|w- \hat{v}\|_{\calL^2(B_r)} \leq \varDelta.
\end{equation}

Finally, to measure the quality of  numerical reconstructions, we introduce the following convenient notation for relative errors:
\begin{equation}\label{def:Err}
		\Err(u, u_0)   :=   \frac{\|u-u_0\|_{\calL^2(B)} }{\|u_0\|_{\calL^2(B)}},
\end{equation}
where $B= B_{\sigma}$ for the case of spatial domain, $B= B_r$ for the case of Fourier domain, and the $\calL^2$-norm is computed using numerical integration.

   \section{Examples}\label{S:Examples}

For our examples, we use the values $\sigma=1$ and  $c = r =10$. We consider the cases $d=1$ and $d=2$.   Our numerical implementations rely on the values 
of $v|_{B_\sigma}$ and $\hat{v}|_{B_r}$ 
on the uniform circumscribed grids of $N^{d}$ points.
The approximate PSWFs $\{\tilde{\psi}_{j,c}\}_{j\geq 0}$  are computed using the software of
\cite{AGD2014}. 
%The approximate eigenvalues $\{\tilde{\mu}_{j,c}\}_{j\geq 0 }$ are taken from the relation   $\tilde{\calF}_c [\tilde{\psi
Then, we find the approximate eigenvalues $\{\tilde{\mu}_{j,c}\}_{j\geq 0}$ 
using the relation $\calF_c [ \psi_{j,c}] = \mu_{j,c}  \psi_{j,c}$.

Table \ref{T:eigen} shows approximations of $\{\tilde{\mu}_{j,c}\}$ for $j$ up to $18$ needed for our computations.
  Note that $n_0 =  \lfloor 2c  / \pi \rfloor =6$ and $\sqrt{\pi/c} \approx 0.560$ and the values $|\tilde{\mu}_{j,c}|$ are in accordance with  
  inequality \eqref{eigen:ineq}.

% 
% Then $\sqrt{\pi/c} \approx 0.560$. 

\begin{table}[h!]
\footnotesize

 \caption{Eigenvalues $\{\mu_{j,c}\}$ for $c=10$ and  $j = 0,\ldots ,18$} \label{T:eigen}
%\centering

\smallskip

\begin{tabular}{c|ccccccccccc}
 $j$  & 0&1&2 &3 & 4&5&6&7 &8& 9 &10
\\
\hline
 $(i)^{-j} \mu_{j,c}$   & 0.793 & 0.793 & 0.793 &  0.792& 0.782& 0.720 &0.526 &0.266   &0.097   &0.029 &0.007
\end{tabular}
 
\begin{tabular}{c|ccccccccc}
 $j$ &11&12&13&14& 15&16&17&18
\\
\hline
 $(i)^{-j} \mu_{j,c}$ & 0.002& $3.7\cdot 10^{-4}$&  $7.1\cdot 10^{-5}$
& $1.3\cdot 10^{-5}$ & $2.2\cdot 10^{-6}$& $3.4\cdot 10^{-7}$& $5.0\cdot 10^{-8}$& $6.9\cdot 10^{-9}$
\end{tabular}

\end{table}

\begin{figure}[h!]
 \centering
 
     \centering
\begin{minipage}[c]{0.9\textwidth}
\centering

 \begin{subfigure}{0.4\textwidth}
 \includegraphics[scale=0.4]{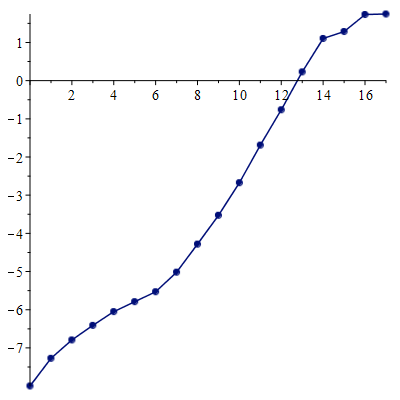}
 \caption{}
\end{subfigure}
\hspace{10mm}
 \begin{subfigure}{0.4\textwidth}
 \includegraphics[scale=0.4]{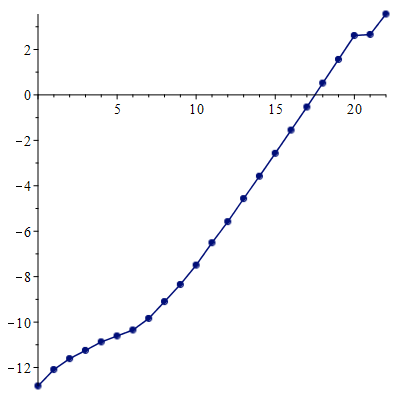}
 \caption{}
\end{subfigure}
\captionsetup{singlelinecheck=off}
 \caption[foo bar]{ The plots of $\log_{10} (\varepsilon_j)$ for (a) $N=129$ and (b) $N=2049$.
%The plot of $\log_10 \eps_j$.
% 
% \begin{itemize}
% 	\item[(a)]  $N = 129$ 
% 	\item[(b)]   $N = 2049$ and  $n=18$.
% 		 The relative errors: $\Err(\tilde{v}_n, v) \approx 0.36$, 
% 	 $\Err(\tilde{v} _{\rm naive}, v) = 0.67$   and  $\Err(\tilde{\calF}[\tilde{v}_n], w) \approx7 \cdot 10^{-9}$,
% 	  $\Err(\tilde{\calF}[\tilde{v} _{\rm naive}], w) \approx  5\cdot 10^{-2}$.
% \end{itemize}
 }
 \label{F:Trust}
 \end{minipage}
\end{figure}

Figure \ref{F:Trust}  shows the corresponding values of $\varepsilon_j$   in the logarithmic scale for $N=129$ and $N=2049$ that we use for our trust criteria \eqref{Trust}.  Clearly,  the numerical calculations involving PSWFs  become unreliable when $n> 12$ for $N=129$
 and when $n> 17$ for $N=2049$.

The following figures  (Figures \ref{1D_nonoise}-- \ref{F:2D_noisy})   illustrate  the reconstruction of  Section \ref{S:2}. 
  The preimages $v$ considered in the present paper
  are rather simple: the sum of characteristic functions of two or three disjoint  objects at distance significantly less than $\pi/r \approx 0.314$.   
 Most importantly,  for all given examples, a proper choice of the regularisation parameter $n$  leads to  super-resolution, that is in this case,  allowing to separate the two objects of the preimage.  
  We also obtain smaller relative errors $\Err$ in $\calL^2$-norm than the naive reconstruction based on formula \eqref{eq:naive} in both  Fourier domain and  spatial domain; see \eqref{def:Err} for the definition of $\Err$.

In this section, we   abbreviate  the notation  $\tilde{\varPhi}_n[w]$ 
  used  in \eqref{def:rmp}  and \eqref{def:Morozov} as  follows:
\begin{equation}\label{num_v}
	\tilde{v}_n := \tilde{\varPhi}_n[w].
\end{equation}
Recall that $w \approx  \calF[v] |_{B_r}$ is the data of Problem \ref{P1} and    $\tilde{\varPhi}_n[\cdot]$  denotes the numerical PSWF reconstruction  via formulas   \eqref{def_vsigma}, \eqref{inverse_1D}, \eqref{inverse_mult}, and \eqref{def:Fnc}.

%Figures 1--Y  show  our numerical examples for the reconstruction of  Section \ref{S:2}. For all these examples, a proper choice of the regularisation parameter $n$  leads to  super-resolution, that is, recovering details of size less than $\pi/r \approx 0.314$. 
%We also obtain a better relative errors in the sense of $\calL^2$-norm than the naive reconstruction based on formula \eqref{eq:naive}.

Figure \ref{1D_nonoise} shows our PSWF reconstructions
$\tilde{v}_n$  with $n = n^*$  defined by  \eqref{def:rmp} from noiseless data $w$, in comparison with preimage $v$ and naive Fourier inversion $\tilde{v} _{\rm naive}$, for $d=1$,  $N=129$, and $N=2049$.  
More precisely, the aforementioned 
data $w$ are noiseless on the uniform grid of $N$
points on $B_r$ in the Fourier domain.
 The most interesting point is that  the reconstruction $\tilde{v}_n$ with $n$ taken according to  the residual minimisation  ($n = n^*$) achieve super-resolution.  
Indeed, the two parts of  $v$ 
are sufficiently distinguished   by $\tilde{v}_n$  even though the 
distance between the two parts   is $\pi/(2r)$.  
The naive reconstruction $\tilde{v} _{\rm naive}$ obscures completely the presence of the two parts.

 Note that even for noiseless data  there still remains "discretisation noise" which gets smaller when $N$ grows. For the example of Figure \ref{1D_nonoise},  it is $0.88\%$ for $N=129$ and $0.06\%$ for  $N=2049$.
Interestingly,  the theoretical suggestion $n = n^*_{\alpha,\delta}$ of \eqref{def:n-star}
 is  close to the residual minimisation choice $n^*$.   Namely,   $n^*_{\alpha,\delta} = 12$  for $N=129$ and  
  $n^*_{\alpha,\delta} = 14$  for $N=129$, where $\alpha=0.75$ and   $\delta$ corresponds  to the   "discretisation noise"    level.

\begin{figure}[h!]
 \centering
 
    \centering
\begin{minipage}[c]{0.9\textwidth}
\centering

 \begin{subfigure}{0.3\textwidth}
 \includegraphics[scale=0.3]{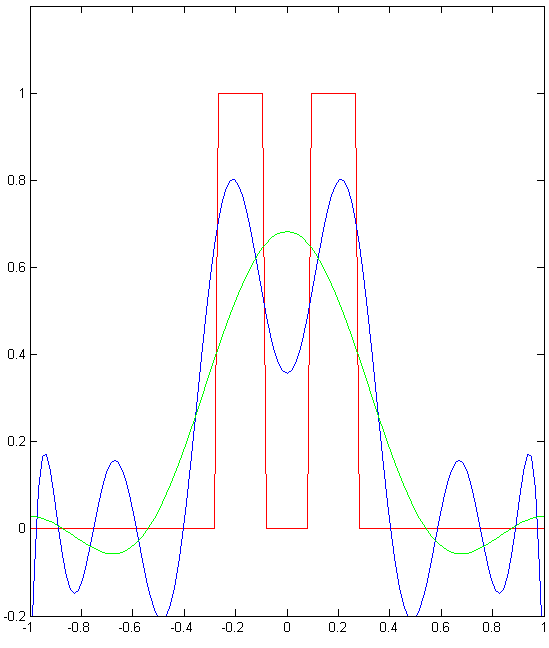}
 \caption{}
\end{subfigure}
\hspace{10mm}
 \begin{subfigure}{0.3\textwidth}
 \includegraphics[scale=0.3]{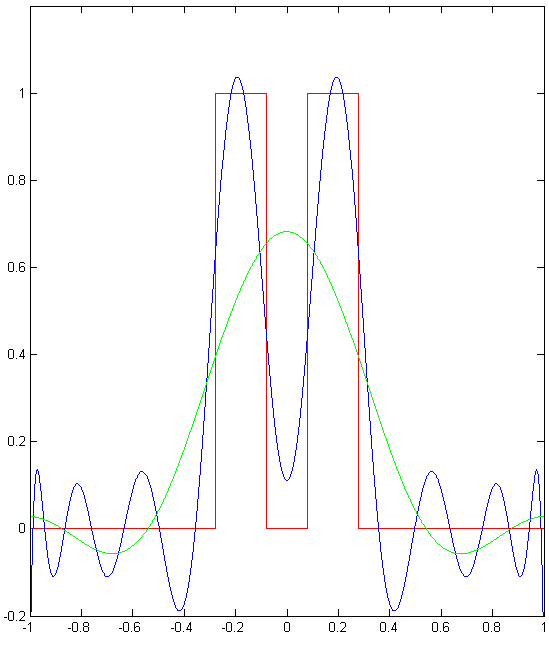}
 \caption{}
\end{subfigure}
\captionsetup{singlelinecheck=off}
 \caption[foo bar]{ \label{1D_nonoise}
 Reconstruction $\tilde{v}_n$(dark blue)   using the residual minimisation from noiseless data $w$, in comparison with preimage $v$(red) and naive Fourier inversion $\tilde{v} _{\rm naive}$(green)  for $d=1$. 
 \begin{itemize}
 	\item[(a)]  $N = 129$ and  $n = n^*=12$.  
 	 The relative errors: $\Err(\tilde{v}_n, v) \approx 0.57$, 
 	 $\Err(\tilde{v} _{\rm naive}, v) \approx 0.71$   and  $\Err(\tilde{\calF}[\tilde{v}_n], w) \approx 4\cdot 10^{-3}$,
 	  $\Err(\tilde{\calF}[\tilde{v} _{\rm naive}], w) \approx 5\cdot 10^{-2}$.
 	     
 	\item[(b)]   $N = 2049$ and  $n= n^*=16$.
 		 The relative errors: $\Err(\tilde{v}_n, v) \approx 0.39$, 
 	 $\Err(\tilde{v} _{\rm naive}, v) \approx 0.67$   and  $\Err(\tilde{\calF}[\tilde{v}_n], w) \approx 4.9 \cdot 10^{-9}$,
 	  $\Err(\tilde{\calF}[\tilde{v} _{\rm naive}], w) \approx  5\cdot 10^{-2}$.
 \end{itemize}
 }

\end{minipage}

\end{figure}
%%%%%%%%%%%%%%%%%%%%%%%%%%%%%%5
%%%%%%%%%%%%%%%%%%%%%%%%%%%%%%%%

%We also observed that  our PSWF reconstructions with $n =n_0$ defined  by \eqref{def:n0}  
%behaves similarly to
%the conventional reconstruction based on formula  \eqref{eq:naive}. 

Figure \ref{1D_nonoise2}  shows our PSWF reconstructions $\tilde{v}_n$
with $n = n_0$ and $n = \tilde{n}_{\varepsilon}+1$, 
where $n_0$ is defined by \eqref{def:n0}  and  
$\tilde{n}_{\varepsilon}$ is defined by \eqref{Trust} with $\varepsilon = 1$,
and    all other parameters are the same as in  Figure  \ref{1D_nonoise}(b). 
Figure \ref{1D_nonoise2}(a) demonstrates the general  phenomenon that 
our reconstruction $\tilde{v}_n$ with $n=n_0$ behaves  similarly to the naive reconstruction $\tilde{v}_{\rm naive}$. Figure \ref{1D_nonoise2}(b) demonstrates that   reconstruction  $\tilde{v}_n$  quickly  becomes  unreliable after  the trust criteria \eqref{Trust} is  violated.

\begin{figure}[h!]
 \centering
   \centering
\begin{minipage}[c]{0.9\textwidth}
\centering

 \begin{subfigure}{0.3\textwidth}
 \includegraphics[scale=0.3]{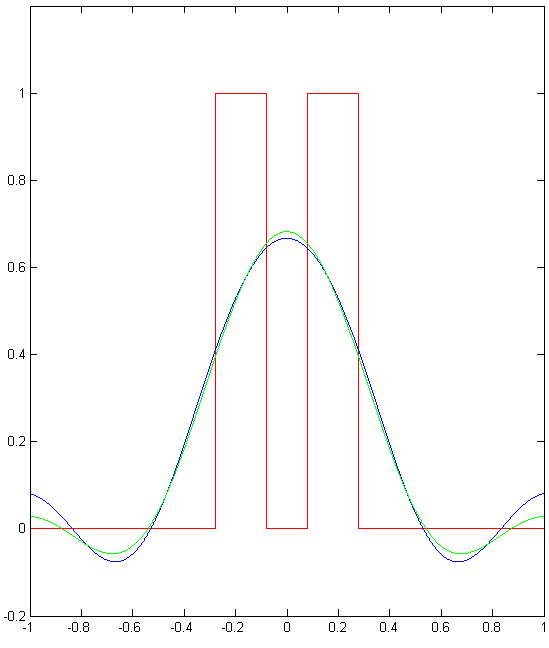}
 \caption{}
\end{subfigure}
\hspace{10mm}
 \begin{subfigure}{0.3\textwidth}
 \includegraphics[scale=0.3]{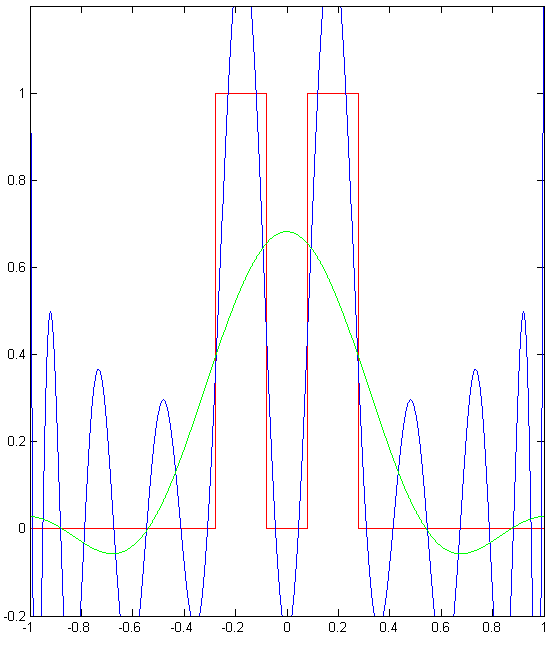}
 \caption{}
\end{subfigure}
\captionsetup{singlelinecheck=off}
 \caption[foo bar]{ \label{1D_nonoise2} 
  Reconstruction $\tilde{v}_n$(dark blue)   from noiseless data  $w$
  in comparison with preimage $v$(red) and naive Fourier inversion $\tilde{v} _{\rm naive}$(green)   for $d=1$ and $N=2049$: 
  (a)~$n=n_0 = 6$,  (b)~$n=\tilde{n}_{\varepsilon}+1 = 18$.
% \begin{itemize}
% 	\item[(a)]  $n=6$  
% 	 The relative errors: $\Err(\tilde{v}_n, v) \approx 0.57$, 
% 	 $\Err(\tilde{v} _{\rm naive}, v) = 0.71$   and  $\Err(\tilde{\calF}[\tilde{v}_n], w) \approx 4\cdot 10^{-3}$,
% 	  $\Err(\tilde{\calF}[\tilde{v} _{\rm naive}], w) \approx 5\cdot 10^{-2}$.
% 	     
% 	\item[(b)]   $N = 2049$ and  $n=16$.
% 		 The relative errors: $\Err(\tilde{v}_n, v) \approx 0.39$, 
% 	 $\Err(\tilde{v} _{\rm naive}, v) = 0.67$   and  $\Err(\tilde{\calF}[\tilde{v}_n], w) \approx 4.9 \cdot 10^{-9}$,
% 	  $\Err(\tilde{\calF}[\tilde{v} _{\rm naive}], w) \approx  5\cdot 10^{-2}$.
% \end{itemize}
 }
\end{minipage}
\end{figure}

%%%%%%%%%%%%%%%%%%%%%%%%%%%%%%5

Figure \ref{F:1Dnoise}  shows our PSWF reconstructions $\tilde{v}_n$ from the noisy data $w \approx \hat{v}|_{B_r}$
with $1.36\%$ of $\calL^2$ random noise
 for $d=1$, $N=129$,  and  $n  \in \left\{ n^*, n^*_{\varDelta}, \frac12 (n^* + n^*_{\varDelta})\right\} $,   where 
  $n^*$, $n^*_{\varDelta}$ are defined by \eqref{def:rmp},  \eqref{def:Morozov} and $\varDelta = 0.0136  \|\hat{v}\|_{\calL^2(B_r)}.$
  In this example, the best reconstruction in the spatial domain  is achieved  when $n =  \frac12 (n^* + n^*_{\varDelta})$. 
   Figure \ref{F:1Dnoise}(c)  illustrates  the well-known fact that   residual minimisation  ($n = n^*$)  may yield explosion
    in the reconstruction from  noisy data.     On the other hand,  
    Morozov's discrepancy principle ($n = n^*_{\varDelta}$)  leads to  a stable reconstruction, but, in our example,  
   it  fails to achieve super-resolution; see    Figure~\ref{F:1Dnoise}(a).

\begin{figure}[h!] 
 \centering
 
  \centering
\begin{minipage}[c]{0.9\textwidth}
\centering

 \begin{subfigure}{0.3\textwidth}
 \includegraphics[scale=0.3]{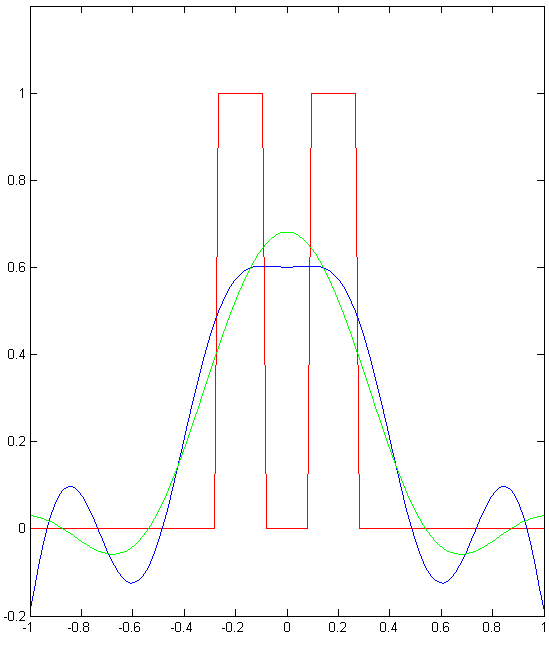}
 \caption{}
\end{subfigure}
\hspace{2mm}
 \begin{subfigure}{0.3\textwidth}
 \includegraphics[scale=0.3]{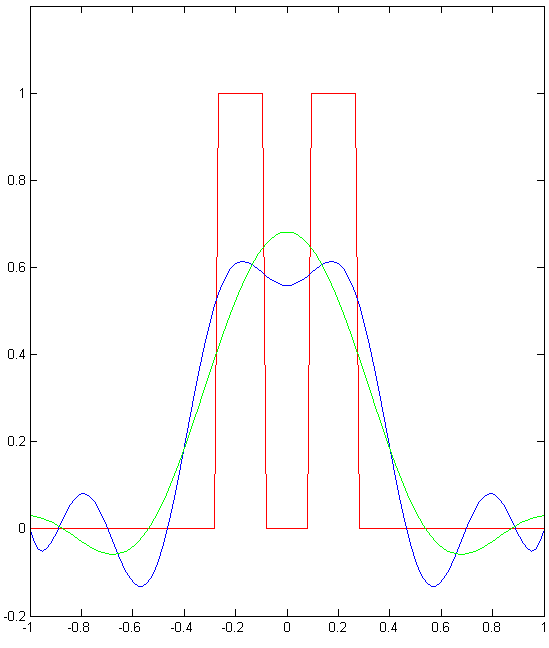}
 \caption{}
\end{subfigure}
\hspace{2mm}
 \begin{subfigure}{0.3\textwidth}
 \includegraphics[scale=0.3]{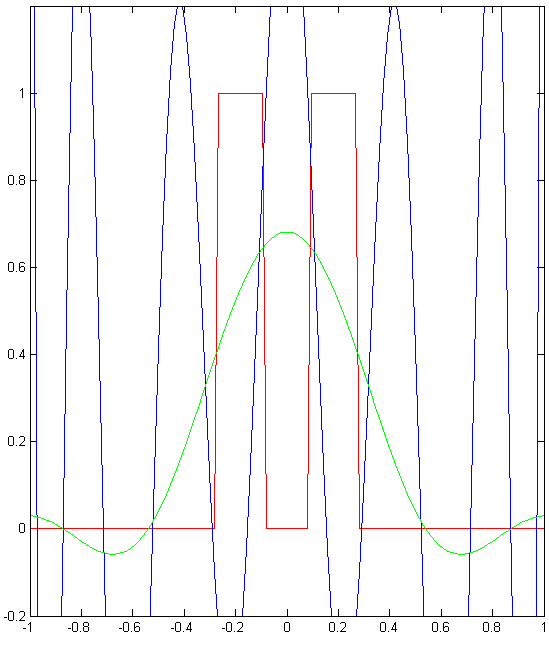}
 \caption{}
\end{subfigure}

\captionsetup{singlelinecheck=off}
 \caption[foo bar]{\label{F:1Dnoise}
 Reconstruction $\tilde{v}_n$ (dark blue)   from noisy data  $w$
 with $1.36\%$ of $\calL^2$ random noise for $d=1$ and $N=129$.
 The preimage $v$  and naive Fourier inversion  $\tilde{v} _{\rm naive}$ are displayed similar to  Figure \ref{1D_nonoise}.
 \begin{itemize}
 	\item[(a)]    The Morozov discrepancy  choice  $n = n_{\varDelta}^*=8$.
 		 The relative errors: $\Err(\tilde{v}_n, v) \approx 0.69$  and  $\Err(\tilde{\calF}[\tilde{v}_n], w) \approx 0.0131$.
 	   
 	   	\item[(b)]    The optimal  choice  $n= \frac12(n_{\varDelta}^* + n^*) = 10$ (best reconstruction in the spatial domain). 
 	   	 The relative errors: $\Err(\tilde{v}_n, v) \approx 0.66$  and  $\Err(\tilde{\calF}[\tilde{v}_n], w) \approx 0.0130$.

 	\item[(c)]  The residual minimisation  choice $n=n^* = 12$.  
 	 The relative errors: $\Err(\tilde{v}_n, v) \approx 2.62$, 
 	  and  $\Err(\tilde{\calF}[\tilde{v}_n], w) \approx 0.0129$.

 \end{itemize}
 }
\end{minipage}
\end{figure}

%%%%%%%%%%%%%%%%%%%%%%%%%
%%%%%%%%%%%%%%%%%%%%%%%%%%%%
%%%%%%%%%%%%%%%%%%%%%%%%%%%
%%%%%%%%%%%%%%%%%%%%%%%%%%
%%%%%%%%%%%%%%%%%%%%%%%%%%%%

Figure  \ref{F:2D_nonoise} and   Figure  \ref{F:2D_nonoise_cross}  show  our PSWF reconstruction 
$\tilde{v}_n$  with $n = n^*$  defined by  \eqref{def:rmp} from noiseless data $w$, in comparison with preimage $v$ and naive Fourier inversion $\tilde{v} _{\rm naive}$ for $d=2$ and  $N=129$.  
In addition,  we implement $\calR^{-1}$ using the filtered back projection (FBP) algorithm
with the angle step of  $2.5^{\circ}$.
Similarly to the one-dimensional example of Figure \ref{1D_nonoise},    the reconstruction $\tilde{v}_n$ with $n$ taken according to  the residual minimisation  ($n = n^*$) achieves super-resolution.
Indeed, the distances  between the three square parts of $v$ are significantly less 
than  $\pi/r$:  two bottom squares are at the distance $0.1$,  while  the top square and any of the bottom squares
are at the distance  $0.05$.
Nevertheless, the three parts of  $v$ 
are sufficiently distinguished   by $\tilde{v}_n$. The naive reconstruction $\tilde{v} _{\rm naive}$ obscures completely the  presence of the three parts.
%
% Note that even for noiseless data  there still remains "discretisation noise" which gets smaller when $N$ grows. For the example of Figure \ref{1D_nonoise},  it is $0.88\%$ for $N=129$ and $0.06\%$ for  $N=2049$.
%Interestingly,  the theoretical suggestion $n = n^*_{\alpha,\delta}$ of \eqref{def:n-star}
% is  close to the residual minimisation choice $n^*$.   Namely,   $n^*_{\alpha,\delta} = 12$  for $N=129$ and  
%  $n^*_{\alpha,\delta} = 14$  for $N=129$, where $\alpha=0.75$ and   $\delta$ corresponds  to the   "discretisation noise"    level.

\begin{figure}[h]
 
 \centering
\begin{minipage}[c]{0.9\textwidth}
\centering
 \includegraphics[scale=0.36]{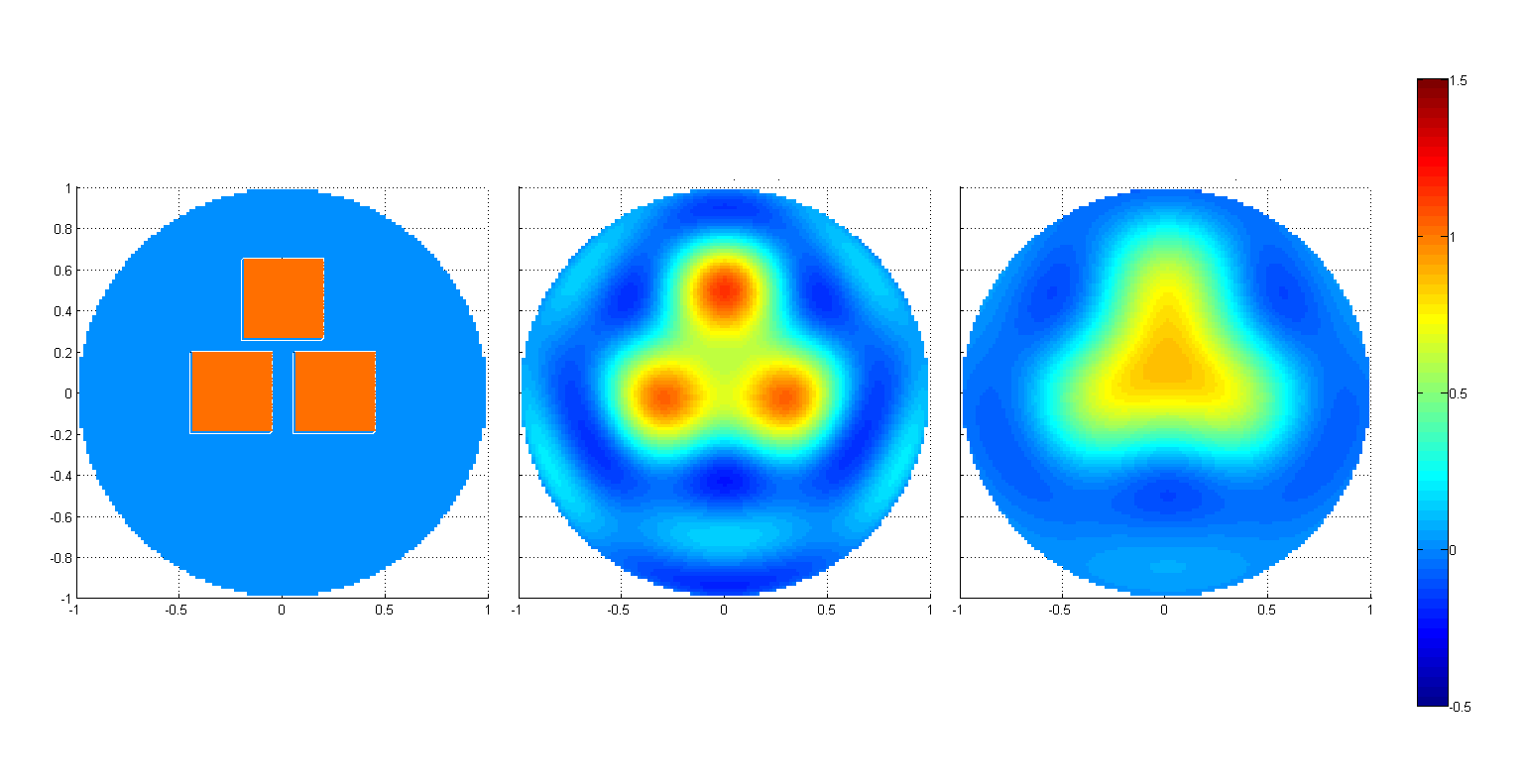}
 
\captionsetup{singlelinecheck=off}
 \caption[foo bar]{\label{F:2D_nonoise}
  Reconstruction  $\tilde{v}_n$(middle) using the residual minimisation from noiseless data $w$  
  in comparison with preimage $v$(left) and  naive Fourier inversion $\tilde{v} _{\rm naive}$(right)
  for $d=2$,  $N=129$, $n= n^* = 9$.
   The relative errors:  $\Err(\tilde{v}_n, v) \approx 0.54$,   $\Err(\tilde{v} _{\rm naive}, v) \approx 0.60$ 
   and    $\Err(\tilde{\calF}[\tilde{v}_n], w) \approx 0.09$,   $\Err(\tilde{\calF}[\tilde{v} _{\rm naive}], w) \approx 0.11$.
%   \begin{itemize} 
% 	\item[(a)]  Preimage $v$.
% 	   
% 	   	\item[(b)]     Reconstruction    	 $\tilde{v}_n$  with $n=8$  using the residual minimisation.   The relative errors: $\Err(\tilde{v}_n, v) \approx 0.54$  and  $\Err(\tilde{\calF}[\tilde{v}_n], w) \approx 0.08$.
% 	     
% 	\item[(c)]     Naive Fourier inversion $\tilde{v} _{\rm naive}$.
% 	The relative errors: $\Err(\tilde{v} _{\rm naive}, v) \approx 0.59$  and  $\Err(\tilde{\calF}[\tilde{v} _{\rm naive}], w) \approx 0.11$.
% \end{itemize}
 }
\end{minipage}
\end{figure}

%%%%%%%%%%%%%%%%%%%%%%%%%%%
%%%%%%%%%%%%%%%%%%%%%%%%%%
%%%%%%%%%%%%%%%%%%%%%%%%%%%%

\begin{figure}[h!]
 \centering

\centering
\begin{minipage}[c]{0.9\textwidth}
\centering

 \begin{subfigure}{0.34\textwidth}
 \includegraphics[scale=0.34 ]{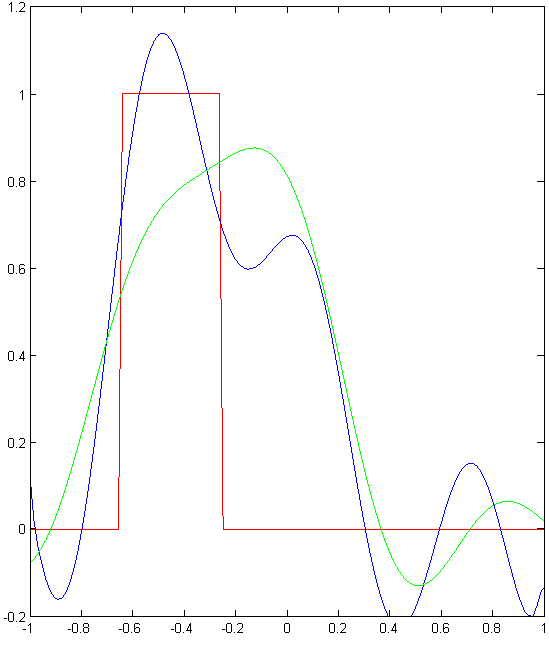}
 \caption{}
\end{subfigure}
 \hspace{10mm}
 \begin{subfigure}{0.34\textwidth}
 \includegraphics[scale=0.34]{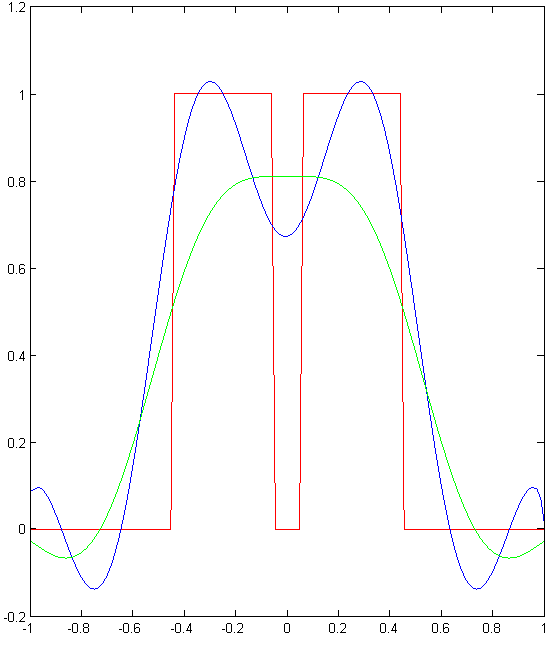}
 \caption{}
\end{subfigure}

\captionsetup{singlelinecheck=off}
 \caption[foo bar]{\label{F:2D_nonoise_cross}
 Cross-sections  of  reconstruction  $\tilde{v}_n$(dark blue)   of Figure \ref{F:2D_nonoise}, in comparison with preimage $v$(red) and naive Fourier inversion $\tilde{v} _{\rm naive}$(green): 
 (a)~along the $y$-axis  ($x=0$);  (b)~along the $x$-axis ($y=0$).
 }
 \end{minipage}
\end{figure}

%%%%%%%%%%%%%%%%%%%%%%%%%%%
%%%%%%%%%%%%%%%%%%%%%%%%%%

   Figure  \ref{F:2D_noisy}  shows  our PSWF reconstruction 
$\tilde{v}_n$  with $n = n^*$  defined by  \eqref{def:rmp} from noisy data $w \approx \hat{v}|_{B_r}$
with $21\%$ of $\calL^2$ random noise  in comparison with preimage $v$ and naive Fourier inversion $\tilde{v} _{\rm naive}$ for $d=2$ and  $N=129$.  In addition, Figure   \ref{F:data} illustrates  the noiseless data $\hat{v}|_{B_r}$ and the noisy data $w$. 
In contrast to  the one-dimensional example of   Figure \ref{F:1Dnoise},    the reconstruction $\tilde{v}_n$ with $n$ taken according to  the residual minimisation  ($n = n^*$)   
works as well as for the noiseless case  shown in Figure  \ref{F:2D_nonoise}.
 Most importantly,  this reconstruction  $\tilde{v}_n$ is rather stable and gives   super-resolution
 even for a considerable level of noise.

%%%%%%%%%%%%%%%%%%%%%%%%%%%
%%%%%%%%%%%%%%%%%%%%%%%%%%
%%%%%%%%%%%%%%%%%%%%%%%%%%%%

\begin{figure}[h]
 \centering

  \centering
\begin{minipage}[c]{0.9\textwidth}
\centering

 \begin{subfigure}{0.4 \textwidth}
 \includegraphics[scale=0.45 ]{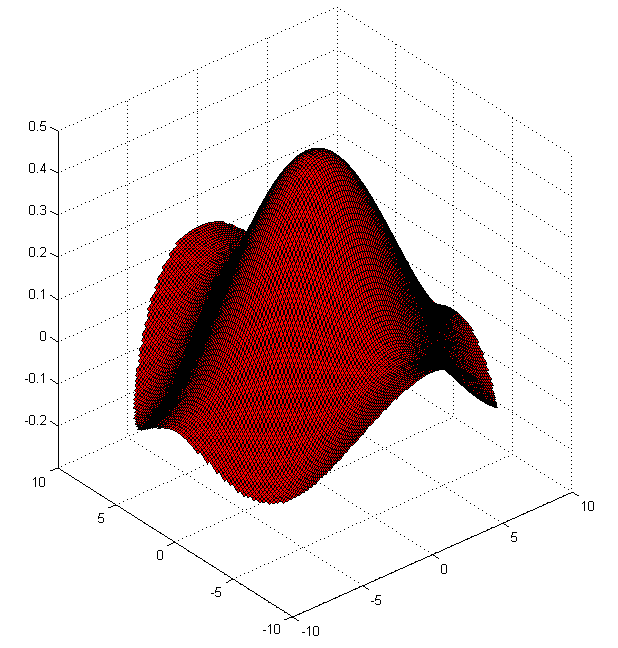}
 \caption{}
\end{subfigure}
\hspace{10mm}
 \begin{subfigure}{0.4  \textwidth}
 \includegraphics[scale=0.45]{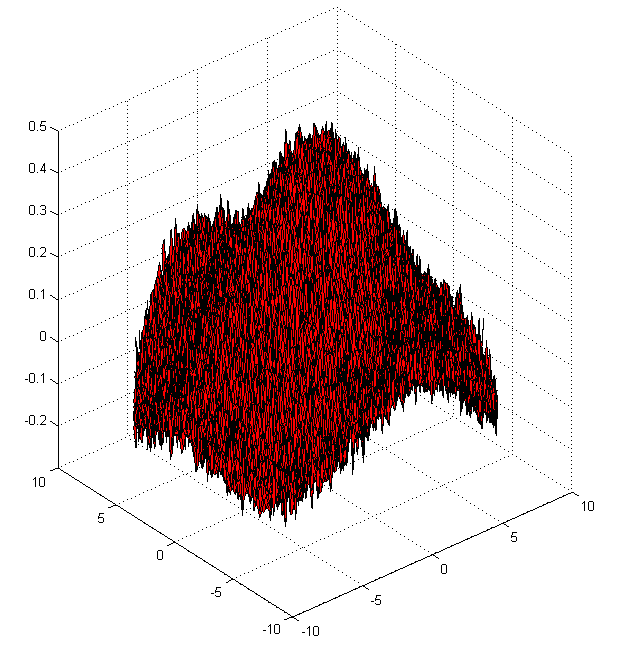}
 \caption{}
\end{subfigure}

\captionsetup{singlelinecheck=off}
 \caption[foo bar]{ \label{F:data}
    The noisy data $w$  with $21\%$  of random $\calL^2$ noise  in comparison with  the  
    noiseless  data $\hat{v}|_{B_r}$
    for the preimage $v$ displayed left on Figure  \ref{F:2D_nonoise} :     (a) the real part of  
    $\hat{v}|_{B_r}$;
    (b) the real part of    $w\approx \hat{v}|_{B_r}$.  
 }
\end{minipage}
\end{figure}

\begin{figure}[h]
 \centering

  \centering
\begin{minipage}[c]{0.9\textwidth}
\centering

 \includegraphics[scale=0.36]{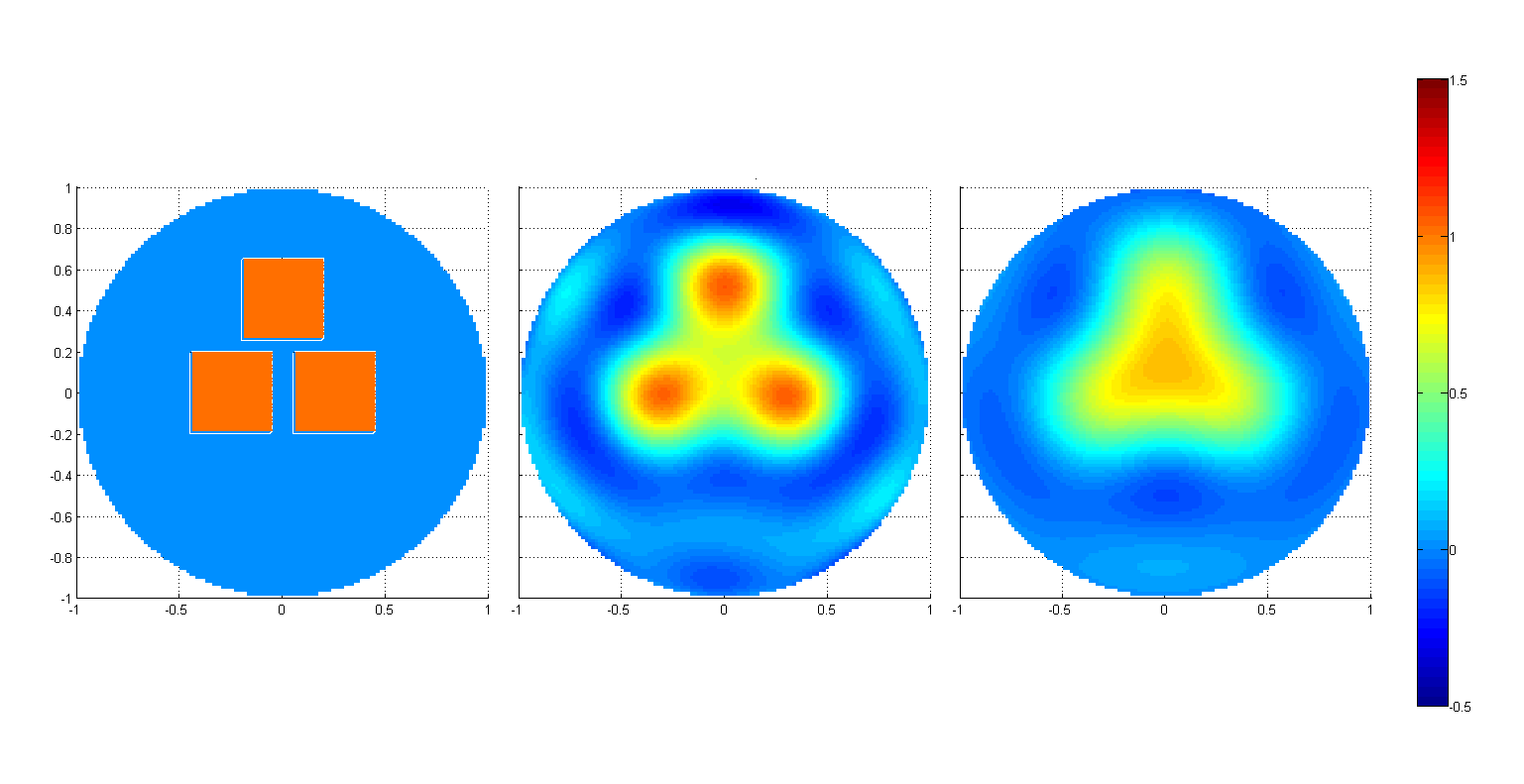}
 
\captionsetup{singlelinecheck=off}
 \caption[foo bar]{ \label{F:2D_noisy}
  Reconstruction  $\tilde{v}_n$(middle) using the residual minimisation from  thenoisy data  displayed 
  on Figure~\ref{F:data}(b)
  in comparison with preimage $v$(left) and  naive Fourier inversion $\tilde{v} _{\rm naive}$(right)
  for $d=2$,  $N=129$, $n= n^* = 8$.
   The relative errors:  $\Err(\tilde{v}_n, v) \approx 0.55$,   $\Err(\tilde{v} _{\rm naive}, v) \approx 0.60$ 
   and    $\Err(\tilde{\calF}[\tilde{v}_n], w) \approx 0.23$,   $\Err(\tilde{\calF}[\tilde{v} _{\rm naive}], w) \approx 0.24$.
%   \begin{itemize} 
% 	\item[(a)]  Preimage $v$.
% 	   
% 	   	\item[(b)]     Reconstruction    	 $\tilde{v}_n$  with $n=8$  using the residual minimisation.   The relative errors: $\Err(\tilde{v}_n, v) \approx 0.54$  and  $\Err(\tilde{\calF}[\tilde{v}_n], w) \approx 0.08$.
% 	     
% 	\item[(c)]     Naive Fourier inversion $\tilde{v} _{\rm naive}$.
% 	The relative errors: $\Err(\tilde{v} _{\rm naive}, v) \approx 0.59$  and  $\Err(\tilde{\calF}[\tilde{v} _{\rm naive}], w) \approx 0.11$.
% \end{itemize}
 }
 \end{minipage}
 \end{figure}

\clearpage

\section{Conclusion}
We implemented numerically  formulas of  \cite{INnotePSWF}  
for finding a compactly supported function $v$ on $\Reals^d$, $d\geq 1$, from its Fourier transform $\calF [v]$ given within the ball $B_r$
(that is, for Problem~\ref{P1}).
Our approach  is based on theoretical and numerical 
results on the prolate spheroidal wave functions, the Radon transform, and regularisation methods.
The present work demonstrates  the numerical efficiency of this approach to Problem~\ref{P1} in its general setting;  including the following  points.

\begin{itemize}
 \item In spite of the exponential instability of the problem, we achieved super-resolution  even for  noisy data by  appropriate choice of the regularisation parameter $n$.  In particular, for $d\geq 2$,  the  approach  works well even for a considerable level of random noise.
 
	\item
	 Our   reconstruction (with appropriate choice of $n$) gives smaller errors  in  $\calL^2$-norm (in both  Fourier domain and spatial domain) than the conventional reconstruction based on formula  \eqref{eq:naive}. 
	
	\item   Our  reconstruction with  	 $n=n_0:= \left\lfloor \dfrac{2 c}{\pi}\right\rfloor$
behaves similarly to
the conventional reconstruction based on formula  \eqref{eq:naive}. 
In our examples, taking $n$ larger than $n_0$ gives better results.
 \end{itemize}

We expect that similar numerical  behaviour  (in particular, super-resolution) 
is also possible for monochromatic inverse scattering (considered, for example, in \cite{ABR2008, BAR2009, HH2001, HW2017, IN2013++, Shurup2022})
and for other  generalisations of Problem \ref{P1}.

%%%%%%%%%%%%%%%%%%%%%%%%%%%%%%%%%%%%%%%%%%%%%%%%%%%%%%%%%%%%%%%%%%%%%%%%%%
%%%%%%%%%%%%%%%%%%%%%%%%%%%%%%%%%%%%%
%%%%%%%%%%%%%%%%%%%%%%%%%%%%%%%%%%%%%
%%%%%%%%%%%%%%%%%%%%%%%%%%%%%%%%%%%%%
%%%%%%%%%%%%%%%%%%%%%%%%%%%%%%%%%%%%%

\section{Acknowledgements}
The work was initiated in the framework of  the internship of  G.V.~Sabinin at  the Centre de Math\'ematique Appliqu\'ees of  Ecole Polytechnique
under the supervision of R.G.~Novikov  in August-October 2021.

\end{document}